\documentclass{amsart}
\usepackage[all, cmtip]{xy}
\usepackage{amssymb}

\usepackage{amsmath}

\usepackage{amsthm}

\usepackage{amscd}

\usepackage{amsfonts}

\usepackage{amssymb}

\newtheorem{prop}{Proposition}[section]

\newtheorem{lemma}[prop]{Lemma}
\newtheorem{cor}[prop]{Corollary}
\newtheorem{conj}[prop]{Conjecture}

\theoremstyle{definition}

\newtheorem{defn}[prop]{Definition}
\newtheorem{eg}[prop]{Example}

\newtheorem{rmk}[prop]{Remark}

\def\p{{\mathbb P}}
\def\cc{{\mathbb C}}

\newcommand\rk{\operatorname{rk}}

\newcommand\crk{\operatorname{crk}}

\newcommand\spec{\operatorname{Spec}}

\begin{document}

\title{Some cases on Strassen additive conjecture}

\author[Youngho Woo]{Youngho Woo}
\address{Department of Mathematical Sciences, KAIST, Taejeon, South Korea}
\email{youngw@kaist.ac.kr}
\begin{abstract} In this article, we verify the additivity for rank of a sum of coprime monomials and bivariate polynomials generalizing the result in (\cite{CCG}). We also show similar results hold for cactus rank.
\end{abstract}

\maketitle
\section{Introduction}
Given a homogeneous polynimial $F \in S^dV$ of degree $d$, the \textit{rank} (or \textit{Waring rank}) $\rk(F)$ of $F$ is the minimal number of linear forms $L_1,...,L_r \in V$ satisfying $F = L_1^d + \cdots + L_r^d$. In \cite{CCG}, the authors show that rank of coprime monomials is the sum of ranks of such monomials. This work verified a particular case of the symmetric version of  Strassen additive conjecture. (See \cite{CCC} and reference therein) 

\begin{conj} If $F \in \cc[x_0, . . . , x_n]$ and $G \in C[y_0, . . . , y_m]$ are homogeneous polynomials of same degree $\geq 2$. Then $\rk(F+G)=\rk(F) +\rk(G)$
\end{conj}

In this paper we verify this conjecture for $F$ and $G$ are sums of several coprime monomials and bivariate polymonimals. (See Corollary \ref{sumofrank} for the precise statement). In the proof, we follows the same track as in \cite{CCG}. This is possible mainly because for a monomial or a bivariate form, the rank attains the lower bound of the inequality (\ref{boundforrank}) in section 2. In fact, we prove the conjecture for those polynomials satisfying this property.
A similar method leads us to the same additivity property for \textit{cactus} rank. (See section 3 and Corollary \ref{crankofsum})

\section{Preliminaries}
\subsection{Apolarity}
Let $V$ be a finite dimensional vector space over the complex field $\cc$ and  
      $S = SV = \bigoplus_{d} S^dV$ be the polynomial algebra over $\cc$ where $S^dV$ is the symmetric $d$-product. Similarly, let $T = SV^* = \bigoplus_{d} S^d(V^*)$ be the polynomial algebra of the dual space $V^*$. Then a polynomial $g \in T$ acts on $S$ as a linear differential operator. Given $F \in S^dV$, the \textit{perp ideal} $F^\perp \subset T$ is the ideal of elements annihilating $F$.  $i.e.$
$$ F^\perp = \{ g \in T | g \cdot F = 0\}.$$
 It is well-known that the ideal $F^\perp$ is a homogeneous ideal such that the quotient $T/F^\perp$ is a zero-dimensional, Gorenstein algebra with socle degree $d$. A subscheme $X \subset \p(V)$ is called \textit{apolar} to $F$   if $I_X \subset F^\perp$ where $I_X$ is the (saturated) homogeneous ideal of the subscheme $X$.  Here we think of $T$ as the homogeneous coordinate ring of $\p(V)$.\\

 Now suppose that $F \in S^dV$. Consider, up to scalar, $d$-th power $L^d$  of a linear form $L \in V$ as an image of the Veronese map $\nu_d : \p(V) \rightarrow \p(S^dV)$ defined by $\nu_d([L]) = [L^d]$. Then the point $[F] \in \p(S^dV)$ is contained in the linear span of the image $\nu_d(\{L_1,...,L_r\})$ if and only if  $F = c_1L_1^d + \cdots + c_rL_r^d$ for some $c_1,...,c_r$. Therefore the rank of $F$ is $$ \min \{r | [F] \in < \nu_d(X) >,  X \subset \p(V) \text{ is reduced, } \dim(X)=0 \text{ and} \deg(X)= r \}$$, where $<\cdot>$ means its linear span in $\p(S^dV)$.\\
The \textit{cactus} rank of $F$ (denoted by $\crk(F)$) is defined similarly,
$$\crk(F) = \min \{r | [F] \in < \nu_d(X) >,\text{ }  \dim(X)=0 \text{ and} \deg(X)= r \}.$$
Note that in this case, we drop the condition "reduced"and consider all $0$-dimensional subschemas $X  \subset \p(V)$ of degree $r$.

\begin{lemma}[Apolarity Lemma \cite{IK}] Let $F \in S^dV$ and $X \subset \p(V)$ be a zero dimensional subscheme. The point $[F] \in \p(S^dV)$ is in the linear span of $\nu_d(X)$ if and only if the subscheme $X$ is apolar to $F$. \textit{i.e.} $I_X \subset F^\perp \subset T$. 
 \end{lemma}

\begin{rmk} Let $F\in S^dV$. According to the apolarity  lemma, we can reformulate the definition of rank and cactus rank 
 \begin{enumerate}
  \item $\rk(F) =\min \{\deg(X) | X\subset \p(V)$ is 0-dimensional, reduced and $I_X \subset F^\perp\} $
  \item $\crk(F) =\min \{\deg(X) | X\subset \p(V)$ is 0-dimensional and $I_X \subset F^\perp\} $
 \end{enumerate}
\end{rmk}
 
\subsection{Degree of a zero-dimensional subscheme} 
Let $X \subset \p(V)$ be a zero-dimensional scheme, and let $I_X \subset T$ be the homogeneous ideal of $X$. Then for a linear form $l \in V^*$ such that $X \cap \{l=0\} = \phi$, we have 
  $$ \deg(X) = \dim_{\cc}T/(I_X+(l))=\sum_{i=0}^{\infty}H(T/(I_X+(l)),i),$$
where $H(T/(I_X+(l)),\cdot)$ is the Hilbert function of $T/(I_X+(l)$.\\
 For a given $F \in S^dV$, if $X$ is apolar to $F$, then we have $\deg(X) \geq \dim_{\cc}T/(F^\perp+(l))$. On the other hand, the ideal $F^\perp +(l)$ defines a zero-dimensional subscheme in $\spec(T)$. The length of this subscheme attains the minimum for a general $l$, by upper semi-continuity. So we conclude that 
\begin{equation}\label{boundforcrank}
 \crk(F) \geq  \dim_{\cc}T/(F^\perp+(l)) \text{ for a \textit{general} linear form } l \in V^*.
\end{equation} 
Now $t \in V^*$ be a linear form, assume that $X$ is reduced. Then the ideal quotient $(I_X : t)$ is the homogeneous ideal of the reduced subscheme $X \setminus \{t=0\}$. As above, if $X$ is apolar to $F$, then we have $\deg(X) \geq \dim_{\cc}T/((F^\perp:t)+(t))$ and
\begin{equation}\label{boundforrank}
 \rk(F) \geq  \dim_{\cc}T/((F^\perp:t)+(t)) \text{ for \textit{any} linear form } t \in V^*.
\end{equation} 
\section{Rank}
\begin{defn} Let $F \in S^d V$ be a homogeneous polynomial of degree $d$, and $t \in V^*$. We will say that the rank of $F$ is \textit{computed by the linear form} $t$ if $$ \rk(F)= \dim_{\cc} T/(F^{\perp} : t) +(t). $$
\end{defn}
\begin{eg}[\textbf{Bivariate case}]\label{rankof2var} Suppose that $\dim V =2$. Then for every $F \in S^d V$, the rank of $F$ is computed by some linear form $ t \in V^*$. In fact, it is well-known that $F^{\perp} = (g_1,g_2)$ for some homogeneous polynomials $g_1,g_2 \in T$ of degree $d_1 \leq d_2$ respectively.  Let $t$ be any linear factor of $g_1$. Then the ideal quotient  $F^{\perp} : t = (g_1/t, g_2$). In fact, it is obvious that $ (g_1/t, g_2) \subset  F^{\perp} : t $. Using an exact sequence
$$ 0 \rightarrow T/(F^{\perp} : t ) \rightarrow T/F^{\perp} \rightarrow T/(F^{\perp} ,t ) \rightarrow 0,$$ we have $\dim T/(F^{\perp} : t ) = (d_1 -1)d_2 = \dim T/(g_1/t, g_2)$ because $t $ is not a factor of $g_2$. Hence $F^{\perp} : t = (g_1/t, g_2)$ holds. Now we have two separate cases
 \begin{enumerate}
 \item ($g_1$ has a double root) Then $\rk(F) = d_2$. Let $t$ be the multiple factor of $g_1$. Since $g_1/t$ still contains $t$,  $\dim_{\cc} T/(F^{\perp} : t) +(t) = d_2$
 \item ($g_1$ has no double roots) Then $\rk(F)=d_1$.  If $d_1 < d_2 $, then take $t$ as one of linear factors of $g_2$ and  $\dim_{\cc} T/(F^{\perp} : t) +(t) = d_1$. If $d_1=d_2$ and $g_2$ has no double roots, then we can choose $g_1' = c_1 g_1 + c_2g_2$ for some constant $c_1,c_2$ such that $g_1'$ has a double root and $F^\perp = (g_1',g_2)$ which belongs to the firtst case. 
 \end{enumerate}

\end{eg}

\begin{eg}[\textbf{Monomial cases}]{\label{rankofmonomials}}(\cite{CCG}, \cite{BBZ}) Let $F = x_0^{a_0}\cdots x_1^{a_n} \in \cc[x_0,...,x_n]$. Let $\{t_0,...,t_n\}$ be the dual basis of $\{x_0,...,x_n\}$. Assume that $a_0 \leq a_i$ for all $i$. Then  $F^\perp = (t_0^{a_0+1},...,t_n^{a_n+1})$ and $\rk(F) = (a_1+1)\cdots (a_n+1)$. Easy to see that $(F^\perp : t_0) +(t_0)=(t_0,t_1^{a_1+1},...,t_n^{a_n+1}) $.  Therefore $\rk(F) =\dim T/(F^\perp : t_0) +(t_0))$ is computed by $t_0$.
\end{eg}

Now we fix some notations. Suppose that we have a decomposition $V = V_1 \oplus V_2 \oplus \cdots \oplus V_m$ for some $m$ and the dual decomposition $V^* = V_1^* \oplus V_2^* \oplus \cdots \oplus V_m^*.$ Then there is an natural idenfication $S(V_i)$ to a subspace of $S=SV$ (and $S(V_i^*)$ to $SV^*$ resp.) for each $i$. By abuse of notation, we will set $S_i = S(V_i)$ and $T_i = S(V_i^*)$. Suppose that $F_i \in S^d(V_i)$ and $t_i \in V_i^*$ for $i=1,...,m$, and let $F = F_1+\cdots+F_m \in S^dV$. Fixing such $F_i$ and $t_i$ for all $i$, we put $J_i =  (F_i^{\perp} : t_i) + ( t_i , V_1^*,...,V_{i-1}^*,V_{i+1}^*,...,V_m^*)$. \textbf{Here $F_i^\perp$ means the perp ideal of $F_i$ \textit{in} $T_i$.} For instance, $J_1$ is the ideal generated by
\begin{itemize}
\item $g\in T_1$ with $t_1g \in F_1^\perp$,
\item $t_1$ and
\item all linear forms in $V_2^*,,...,V_m^*$.
\end{itemize} under the inclusion $T_i \subset T $. \\

\begin{lemma}Notation as above.
 $(F^\perp:(t_1,...,t_m)) \subset J_1\cap \cdots \cap J_m.$
\end{lemma}
PROOF. By symmetry, it is enough to show that $ F^\perp :t_1 \subset  (F_1^{\perp} : t_1) + ( t_1 , V_2^*,...,V_m^*).$ Let $h \in (F^\perp : t_1)$. We can write $h = h_1 + h_2$ where $h_1 \in T_1 $ and $ h_2 \in (V_2^*,...,V_m^*)$. We claim that $h_1 \in (F_1^{\perp} : t_1) $ . By definition, $t_1h \in F^\perp$. Since  $t_1h_2 \in F^\perp$, we have $t_1h_1 = t_1h -  t_1h_2 \in F^\perp \cap T_1 = F_1^\perp $ and hence $h_1 \in ( F_1^\perp : t_1)$. The inclusion follows.\\

\begin{prop}Suppose that we have a decompsition $V = V_1 \oplus V_2 \oplus \cdots \oplus V_m$ for some $m$. Let $F_i \in S^d V_i$ for $1 \leq i \leq m$ and for $d\geq2$, let $F  = F_1 + ... + F_m \in S^d V$. If the rank of $F_i$ is computed by some $t_i \in V_i^*$  for all $i$, then  $\rk(F)=\rk(F_1) + ...+ \rk(F_m)$.
 \end{prop}
 
 The proof proceeds similarly as in \cite{CCG}\\
 PROOF. We will prove that 
 \begin{equation}\label{i1}\rk(F) \geq \rk(F_1) +\cdots+ rk(F_m)
 \end{equation} because the opposite inequality is obvious.\\ 
Let $r = \rk(F)$. Then by apolarity lemma, there is a finite set $X \subset \p(V)$ of degree $r$ such that the (saturated) homogeneous ideal $I_X \subset F^{\perp}$. Then the quotient ideal  $ I_X : (t_1,t_2,..,t_m)$ is the homogeneous ideal of the subset $X'$ of $X$ not lying on $\{ t_1=t_2=\cdots=t_m=0\}$. Then a linear form $\lambda_1t_1 +\cdots \lambda_m t_m \in V^*$ for general $\lambda_1,...,\lambda_m \in \cc$ does not vanish any point of $X'$. By replacing $t_i$ by $\lambda_it_i$, we may assume that $ X' \cap \{ t_1+\cdots+ t_m=0\} = \phi$. 
By the previous  lemma and the fact $(t_1 + \cdots + t_m) \subset   J_1\cap \cdots \cap J_m$, we have 
 \begin{equation}\label{i2}I_{X'} + (t_1 + \cdots + t_m) \subset   J_1\cap \cdots \cap J_m
  \end{equation}\\
Now we prove two claims.\\

\textbf{Claim 1.} $ \dim_{\cc} T/(J_1 \cap \cdots \cap J_m) = \rk(F_1) + \cdots \rk(F_m) -(m-1).$\\
For each $i =1,...,2$, one has $T/J_i \cong T_i /((F_i^\perp : t_i) + (t_i))$, and hence $ \dim T/J_i = \rk(F_i)$. From the exact sequence
$$ 0 \rightarrow T/(J_1\cap \cdots \cap J_{m}) \rightarrow  T/(J_1\cap\cdots \cap J_{m-1}) \oplus T/J_m \rightarrow T/(J_1\cap \cdots \cap J_{m-1}+ J_m)\rightarrow 0,$$
and $T/(J_1\cap\cdots\cap J_{m-1}+ J_m) \cong \cc$, the assertion follows by induction.\\

\textbf{Claim 2.}  $t_1,...,t_m$ are linearly independent modulo $I_{X'}$.\\
Suppose $\mu_it_1 + \cdots + \mu_mt_m \in I_{X'}$ for some $\mu_1,...,\mu_m \in \cc$. Since this already vanishes along $\{t_1=\cdots=t_m =0 \}$, we have $\mu_1t_1 + \cdots \mu_m t_m \in I_X \subset F^\perp.$ Then $(\mu_1t_1 + \cdots \mu_m t_m)\cdot F =  \mu_1t_1\cdot F_1 + \cdots + \mu_mt_m\cdot F_m =0$. Note that each polynomial $t_i\cdot F_i \in S_i(=S(V_i))$ has degree $\geq 1$ or $t_i\cdot F_i=0$  for $i=1,...,m$. Since $t_i $ compute the rank of $F_i$, $t_i\cdot F_i \neq 0$. (Otherwise $1\in (F_i:t_i)$.) Thus $\mu_1=\cdots=\mu_m =0$\\

From the claim 2 and the inclusion (\ref{i2}), we conclude that
\begin{equation}\label{i3}
H(T/(I_{X'} + (t_1 +\cdots+ t_m)),1) \geq H(T/(J_1 \cap \cdots \cap J_m, 1)+ (m-1)
\end{equation}
 since $t_1,....,t_m \in J_1\cap\cdots \cap J_m$. \\
Now we verify the inequality (\ref{i1}).

 \begin{align*}
r \geq & \deg(X' )  =   \dim T/((I_{X'} + (t_1 +\cdots+ t_m))  \\
=        & \sum_{j = 0}^{\infty} H(T/(I_{X'} + (t_1 +\cdots+ t_m)),j) \\
=        &H(T/(I_{X'} + (t_1 + \cdots + t_m)),1)  +  \sum_{j \neq 1}H(T/(I_{X'} + (t_1 +\cdots +t_m)),j)\\
\geq   &H(T/(J_1 \cap\cdots\cap J_m, 1)+(m-1) +  \sum_{j \neq 1}H(T/(J_1 \cap\cdots\cap J_m, j) \text{ by (\ref{i2}) and (\ref{i3})} \\
=        &\dim_{\cc} T/(J_1 \cap \cdots \cap J_m) +(m-1) \\
=        &\rk(F_1) + \cdots \rk(F_m)  \text{ by the claim 1} \\
 \end{align*}

\begin{cor}\label{sumofrank} Let $F_i \in S^dV_i$ $(d \geq 2)$ be a monomial or a bi-variate polynomial for each $i$. Then $\rk(F_1 + \cdots + F_m) = \rk(F_1) + \cdots +\rk(F_m)$.
\end{cor}

\section{Cactus Rank}

\begin{defn} Let $F \in S^d V$ be a homogeneous polynomial of degree $d$. We will say that the cactus rank of $F$ is \textit{computed by a general linear form}  if $$ \crk(F)= \dim_{\cc} T/(F^{\perp}  + (l)). $$ for a general linear form $l \in V^*$.
\end{defn}

\begin{eg} Suppose that $\dim V =2$. Then for every $F \in S^d V$, the cactus rank of $F$ is computed by a general linear form. As in the example (\ref{rankof2var}), $F^{\perp} = (g_1,g_2)$ for some homogeneous polynomials $g_1,g_2 \in T$ of degree $d_1 \leq d_2$ respectively. Then $\crk(F) = d_1 = \dim_{\cc} T/((g_1,g_2)  + (l)) $ for a general  $l \in V^*$ .\\
\end{eg} 

\begin{eg}{\label{crankofmonomials}}(\cite{RS}) As in the example (\ref{rankofmonomials}), let $F = x_0^{a_0}\cdots x_1^{a_n} \in \cc[x_0,...,x_n]$. Here we give an additional assumption that $a_0+\cdots+a_{n-1} \leq a_n$. Then  $F^\perp = (t_0^{a_0+1},...,t_n^{a_n+1})$ and $\crk(F) = (a_0+1)\cdots (a_{n-1}+1)$ is computed by a general linear form. In fact,  for general $ \lambda_0,...,\lambda_n  \in \cc $
\begin{align*} 
\dim_{\cc} &T/(F^\perp + (\lambda_0t_0+\cdots \lambda_n t_n ))\\
=&\dim_{\cc} \cc[x_0,...,x_{n-1}]/(t_0^{a_0+1},...,t_{n-1}^{a_{n-1}+1} ,(t_0+\cdots+ t_{n-1})^{a_n +1}) \\
=&\dim_{\cc} \cc[x_0,...,x_{n-1}]/(t_0^{a_0+1},...,t_{n-1}^{a_{n-1}+1})\\
=& (a_0+1)\cdots (a_{n-1}+1)
\end{align*}
because $(t_0+\cdots+ t_{n-1})^{a_n +1} \in (t_0^{a_0+1},...,t_{n-1}^{a_{n-1}+1})$ provided that  $a_0+\cdots+a_{n-1} \leq a_n$
\end{eg}


From now on, fix $F \in S^dV$ for $d \geq 2$. Suppose that $ F = F_1 + F_2$ for some $F_1 \in S^dV_1$ and $F_2 \in S^dV_2$ under the decomposition $V = V_1 \oplus V_2 $.
\begin{lemma}Let $l_i \in V_i^* $ for $i=1,2$. Suppose  that $d\geq2$ and  $l_i \notin F_i ^\perp$. Then 
 \begin{enumerate}
  \item $F^\perp + (l_1 +  l_2) \subsetneq [(F_1^\perp + (l_1) + (V^*_2 )]  \cap [F_2^\perp + (l_2) + (V_1^*)]$
  \item $\dim T/ (F^\perp + (l_1 + l_2)) \geq \dim T_1/(F_1^\perp + (l_1))+ \dim T_2/(F_2^\perp + (l_2)) $
 \end{enumerate}
\end{lemma}

PROOF. At first, we show that $F^\perp \subset F_1^\perp + (l_1) + (V_2^*).$ Let $h \in F^\perp$ be a nonzero homogenous polynomial. Since the set $F_1^\perp + (V_2^*)$ contains all monomials of degree $> d$, we may assumt that $\deg(h) \leq d$. Now we can write $h = h_1 + h_2$ for $ h_1 \in T_1$ and $h_2 \in V_2^*$. Since $h\cdot F=0$, we have $h_1\cdot F_1 + h_2\cdot F_2 = 0$. If $\deg(h) < d$, then $h_1\cdot F_1= h_2\cdot F_2 = 0$ and hence  $h \in F_1^\perp + (V_2^*)$. On the other hand, if $\deg(h) = d$, then   both $h_1\cdot F_1$ and  $h_2\cdot F_2$ are constant. Since $l_1 \notin F_1^\perp$, there is $g \in T_1$ of degree $d-1$ with $g\cdot (l_1\cdot F_1)$ is a nonzero constant. Then $h_1 - \frac{(h_1\cdot F)}{gl_1\cdot F_1}gl_1 \in F_1^\perp$. So we have $ h = (h_1 - \frac{(h_1\cdot F)}{gl_1\cdot F_1}gl_1) + \frac{(h_1\cdot F)}{gl_1 \cdot F_1}gl_1  + h_2 \in F_1^\perp + (l_1) + (V_2^*)$. Up to symmetry, the inclusion of the lemma follows. This inclusion must be strict because $l_1\notin F^\perp + (l_1 +  l_2) $. Otherwise, there is a constant $\lambda$ such that $ (1-\lambda)l_1 - \lambda l_2 \in F^\perp$, and hence $ (1-\lambda)l_1 \cdot F_1 - \lambda l_2\cdot F_2 =0$. But this is impossible since $\deg(F) \geq 2$ and $l_i \notin F_i^\perp $ for $i=1,2$. This proves (1).\\
Now Put $J_1 = F_1^\perp + (l_1) + (V^*_2)$, and $J_2 = F_2^\perp + (l_2) + (V_1^*)$ as in the rank case. By (1), we have 
\begin{align*}
   & \dim T/(F^\perp + (l_1 +  l_2))\\
> & \dim T/J_1\cap J_2 \\
= & \dim T/J_1 + \dim T/ J_2 - \dim T/(J_1+J_2) \\
=  & \dim T_1/(F_1^\perp + (l_1)) + T_2/(F_2^\perp + (l_2)) - 1\\
\end{align*}

\begin{prop}Suppose that the cactus rank of $F_i$ is computed by a general linear form for $i =1,2$. Then
 $\crk(F) =\crk(F_1) + \crk (F_2)$ and this is computed by a general linear form.
\end{prop}
PROOF. We always have $\crk(F) \leq \crk(F_1) + \crk (F_2)$. By assumption, $ \crk(F_i) = \dim T_i /(F_i^\perp + (l_i))$ for a general linear form $l_i$ for $i=1,2$.  Then $l_i \notin F_i ^\perp$.  By apolarity lemma, there is a zero dimensional  subscheme $X \subset \p(V)$ such that  $ I_X \subset F^\perp$ and $\deg(X)=\crk(F)$. Since $l_1+l_2$  is a general linear form in $V^*$, we have
\begin{align*}
\deg(X) = & \dim T/(I_X + (l_1 + l_2)) \\
\geq  & \dim T/ (F^\perp + (l_1+l_2)) \\
\geq  & \dim T_1/(F_1^\perp + (l_1))+ \dim T_2/(F_2^\perp + (l_2)) &\text{(by the lemma)} \\
 =      &\crk(F_1) + \crk(F_2)\\
\end{align*}

\begin{cor}\label{crankofsum}Let $F_i \in S^dV_i$ $(d \geq 2)$  be a monomial satisfying the assumption in example \ref{crankofmonomials} or a bi-variate polynomial for each $i=1,...,m$. Then $\crk(F_1 + \cdots + F_m) = \crk(F_1) + \cdots +\crk(F_m)$.
\end{cor}

\end{document}